\documentclass[a4paper,11pt,oneside]{article}
\usepackage[english]{babel}

\usepackage[dvips]{graphicx}
\usepackage[svgnames]{xcolor}
\colorlet{backside}{Bisque!90!LightGray}
\colorlet{foldline}{Red}
\colorlet{vector}{RoyalBlue!100!black}
\colorlet{halfgray}{DarkGrey}

\usepackage{amssymb}
\usepackage[tbtags]{amsmath}
\usepackage{amsthm}
\usepackage{icomma}

\theoremstyle{definition}
\newtheorem{definition}{Definition} 
\newtheorem{incidence}{Incidence}

\newtheorem{operation}{Operation}

\theoremstyle{plain}
 
\newtheorem{theorem}{Theorem} 

\usepackage[T1]{fontenc}
\usepackage[osf,sc]{mathpazo}
\usepackage[euler-digits,euler-hat-accent]{eulervm}

\usepackage[top=32mm]{geometry}

\usepackage{caption}
\usepackage{booktabs}
\usepackage[flushleft]{threeparttable}

\usepackage{tabularx}

\usepackage[authoryear,round,longnamesfirst]{natbib}
\usepackage{fancyhdr}
\usepackage{hyperref}
\hypersetup{
  breaklinks = true,
  linktocpage=true,
  linkcolor  = MidnightBlue,
  citecolor  = MidnightBlue,
  urlcolor   = MidnightBlue,
  colorlinks = true,
}
\usepackage{breakcites}    

\usepackage{pst-all} 
\newpsobject{showgrid}{psgrid}{subgriddiv=1,griddots=10,gridlabels=6pt}
\psset{arrowscale=2}

\title{On the elementary single-fold operations of origami: reflections and incidence constraints on the plane}
\author{Jorge C. Lucero\thanks{Dept.\ Computer Science, University of Bras\'{i}lia, Brazil. E-mail: \href{mailto:lucero@unb.br}{lucero@unb.br} }} 
\date{\today}

\graphicspath{{./FG/}}

\begin{document}

\maketitle

\begin{abstract} 
\noindent 
This article reviews the so-called ``axioms'' of origami (paper folding), which are elementary single-fold operations to achieve incidences between points and lines in a sheet of paper. The geometry of reflections is applied, and exhaustive analysis of all possible incidences reveals a set of eight elementary operations. The set includes the previously known seven ``axioms'', plus the operation of folding along a given line. This operation has been ignored in past studies because it does not create a new line. However, completeness of the set and its regular application in practical origami dictate its inclusion. Formal definitions and conditions of existence of solutions are given for all the operations.    
\end{abstract}

\section{Introduction}

Three decades ago, \citet{Justin1986} introduced a set of 7 elementary single-fold operations which have become known as the ``axioms'' of origami, the Japanese art of paper folding \citep{Alperin2000,Alperin2006,Ghourabi2013}. Each operation is defined by one or more alignments (incidences) between points and lines on a sheet of paper, that must be achieved with a single fold.  The number of solutions of each operation must be finite; however, depending on the relative position of the given points and lines, some of the operations may have none, one or multiple solutions. It has been shown that the set of operations constitutes a more powerful geometrical tool than the combination of straight edge and compass \citep{Alperin2000}. For example, the operations allow for the trisection of arbitrary angles \citep{Hull1996}, solving the problem of duplicating  the cube \citep{Messer1986}, constructing heptagons \citep{Geretschlaeger1997} and solving cubic and fourth order equations \citep{Alperin2000,Geretschlager1995}, all of which may be not done by straight edge and compass alone. 

Justin's work \citep{Justin1986} seems to have been overlooked at its time, and the same fold operations have been rediscovered later and expressed under various forms by \citet{Huzita1989}, Hatori \citep[2001, according to][]{Alperin2006} and other enthusiasts of origami mathematics \citep[e.g.,]{Alperin2000,Auckly1995, Geretschlager1995,Kawasaki2011,Martin1998}. The operations are popularly known today as Huzita's axioms, Huzita-Justin's axioms or Huzita-Hatori's axioms.\footnote{In his work, \citet{Huzita1989} listed six of Justin's fold operations. The seventh was rediscovered by Hatori, in 2001 (according to \citep{Alperin2006})} Let us note that the designation as ``axioms'' is not correct since some of the operations may be derived from others, and further, some of them may not be possible depending on the configuration of given points and lines. 
It has been claimed that the set is complete, in the sense that it contemplates all possible alignments between points and lines with a finite number of solutions, and excluding redundant alignments \citep{Alperin2006}. 

More recently, \citet{Kasem2011} showed several inconsistencies in the operations owing to the rather imprecise form in which they had been stated. The inconsistencies included impossibility of some folds, infinite solutions and superfluous conditions. Seeking a more rigorous treatment, \citet{Ghourabi2013} expressed the operations in formal algebraic terms and analyzed their number of solutions and conditions of existence. Such a formalization is a necessary step for adapting folding techniques to industrial applications. In fact, recent years have seen a surge of applications of origami to science and technology, e.g., in aerospace and automotive technology \citep{Cipra2001}, materials science \citep{Zhou2016}, biology \citep{Mahadevan2005}, civil engineering \citep{Filipov2015}, robotics \citep{Felton2014} and acoustics \citep{Harne2016}. Further, a number of computational systems of origami simulation have been developed \citep{Ida2006a}. 
   
The present article follows the call for normalization and analyzes the elementary operations by applying the geometry of reflections. Folding a sheet of paper along a straight line superposes the paper on one side of the fold line to the other side. 
As a result, all points and lines on each side of fold line are reflected across it  onto the
other side \citep{Martin1998}. In fact, the geometry of paper folding may be reproduced by using a semi-reflective mirror called ``Mira'' \citep{Demaine2007}. A full geometrical characterization of Mira constructions has been  given by \citet{Emert1994} in terms of ``primitive actions'', which are equivalent to the fold operations studied here.

First, the analysis will determine all possible incidences between given points and lines on a plane that may be achieved by a reflection. Next, it will derive all possible fold operations that may be defined so as to satisfy combinations those incidences. In this way, a total of eight elementary operations will be obtained, i.e., one more than the previous set, where the new operation is to fold along a given line. This operation has been commonly disregarded by previous studies on the argument that it does not create a new line; however, completeness of the set demands its inclusion. Further, it has applications in actual origami folding, as will be discussed later (Section \ref{foldingalong}).

\section{Reflections on a plane}
\label{reflections}

The medium on which all folds are performed is assumed to be an infinite Euclidean plane \citep{Geretschlager1995}. Points are denoted by capital letters ($P$, $Q$ etc.), lines by small letters ($m$, $n$ etc.), except the fold line which is denoted by the special symbol $\chi$, and $P\in m$ means that point $P$ is on line $m$. 

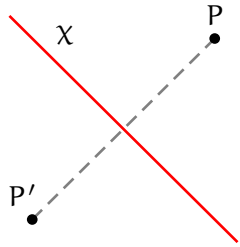
\begin{figure}[!htb]
\centering
\begin{pspicture}(1,1.5)(5,5)
%\showgrid
\psline[linecolor=foldline,linewidth=1pt](1.5,4.5)(4.5,1.5)
\psline[linestyle=dashed,linecolor=gray,linewidth=1pt](1.8,1.8)(4.2,4.2)
\qdisk(1.8,1.8){2pt}
\uput[110](1.8,1.8){$P'$}
\qdisk(4.2,4.2){2pt}
\uput[90](4.2,4.2){$P$}
\uput[45](2,4){$\chi$}
\end{pspicture}
\caption{Reflection of point $P$ in line $\chi$.} 
\label{reflx}
\end{figure}

A reflection is defined as follows \citep{Martin1998}:

\begin{definition}\label{defreflx} Given a line $\chi$, the reflection $\mathcal{F}$ in $\chi$ is the mapping on the set of points in the plane such that for point $P$
\[
\mathcal{F}(P)=
\left\{
\begin{array}{ll}
P & \text{if $P\in \chi$},\\
P' & \text{if $P\notin \chi$ and $\chi$ is the perpendicular bisector of segment $\overline{PP'}$}
\end{array}
\right.
\]
(see Fig.~\ref{reflx}).
\end{definition}

It is easy to see that $\mathcal{F}(P)=P'$ iff $\mathcal{F}(P')=P$.

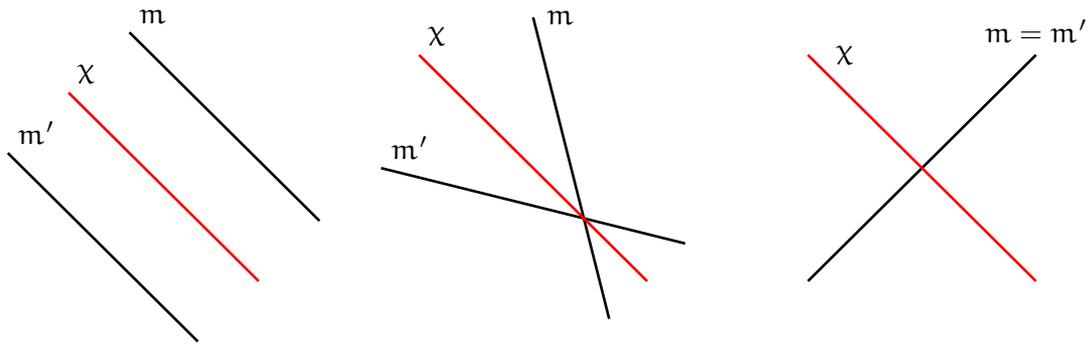
\begin{figure}[!htb]
\begin{pspicture}(1,0.5)(6,5.5)
%\showgrid
\psline[linecolor=foldline,linewidth=1pt](2,4)(4.5,1.5)
\psline[linewidth=1pt](5.3,2.3)(2.8,4.8)
\psline[linewidth=1pt](1.2,3.2)(3.7,0.7)
\uput[45](2.8,4.8){$m$}
\uput[45](1.2,3.2){$m'$}
\uput[45](2,4){$\chi$}
\end{pspicture}
\begin{pspicture}(1,0.5)(6,5.5)
%\showgrid
\psline[linewidth=1pt](3,5)(4,1)
\psline[linewidth=1pt](1,3)(5,2)
\psline[linecolor=foldline,linewidth=1pt](1.5,4.5)(4.5,1.5)
\uput[0](3,5){$m$}
\uput[45](1,3){$m'$}
\uput[45](1.5,4.5){$\chi$}
\end{pspicture}
\begin{pspicture}(1,0.5)(5,5.5)
%\showgrid
\psline[linewidth=1pt](1.5,1.5)(4.5,4.5)
\psline[linecolor=foldline,linewidth=1pt](1.5,4.5)(4.5,1.5)
\uput[90](4.5,4.5){$m=m'$}
\uput[45](1.75,4.25){$\chi$}
\end{pspicture}
\caption{Reflection of line $m$ in line $\chi$. Left: $m\parallel \chi$. Center: $m\nparallel\chi$. Right: $m\bot\chi$} 
\label{linereflx}
\end{figure}

Reflection of a line $m$ in $\chi$ is obtained by reflecting every point in $m$. Therefore, $\mathcal{F}(m)=\{\mathcal{F}(P)|\,P\in m\}$.
Let $m'=\mathcal{F}(m)$ and  consider the following cases:

\begin{enumerate}

\item If $m$ and $\chi$ are parallel ($m\parallel \chi$), then $m\parallel m'$ (Fig.~\ref{linereflx}, left).

\item\label{it2} If $m$ and $\chi$ are not parallel ($m\nparallel \chi$)  then line $\chi$ is a bisector of the angle between $m$ and $m'$ (Fig.~\ref{linereflx}, center).

\item If $m=\chi$, then every point in $m$ is its own reflection, and therefore $m'=m$.

\item \label{R2.4} If $m$ and $\chi$ are perpendicular ($m\bot\chi$), then the reflection of every point of $m$ is also on $m$; therefore, $m=m'$ (Fig.~\ref{linereflx}, right). Note also that $\chi$ divides $m$ into two halves, and each half is reflected onto the other. Thus, for every point $P$ on $m$ and not on the intersection with $\chi$, $\mathcal{F}(P)\neq P$. 

\end{enumerate}

\section{Incidence constraints}

Elementary single-fold operations are defined in terms of incidence constraints  between pairs of objects (points or lines) that must be satisfied with a fold \citep{Alperin2006, Ghourabi2013,Justin1986}. Each constraint involves an object $\alpha$ and the image $\mathcal{F}(\beta)$ of an object $\beta$ (including the case $\alpha=\beta$) by reflection in the fold line.  The symmetry of the reflection mapping implies that all incidence relations are also symmetric.

A total of six different incidences are possible on a plane, an they are defined and analyzed in the next subsections (see also Table \ref{table0}). In order to facilitate the posterior definitions of the fold operations, incidences involving distinct objects (i.e., $\alpha\neq \beta$) are distinguished from those involving the same object (i.e., $\alpha= \beta$).

\subsection{Incidences involving distinct objects}

\begin{incidence}
$\mathcal{F}(P)=Q$, with $P\neq Q$.
\label{I1}
\end{incidence}

\noindent In this incidence, the reflection of a given point $P$ coincides with another given point $Q$. According to Definition \ref{defreflx}, its solution is the unique fold line $\chi$ which is the perpendicular bisector of segment $\overline{PQ}$. 

\begin{incidence}
$\mathcal{F}(m)= n$, with $m\neq n$.
\label{I2}
\end{incidence}

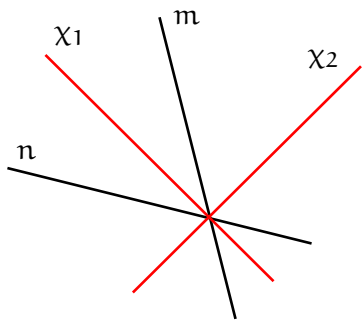
\begin{figure}[!htb]
\centering
\begin{pspicture}(1,0.5)(6,5.5)
%\showgrid
\psline[linewidth=1pt](3,5)(4,1)
\psline[linewidth=1pt](1,3)(5,2)
\psline[linecolor=foldline,linewidth=1pt](1.5,4.5)(4.5,1.5)
\psline[linecolor=foldline,linewidth=1pt](2.65,1.35)(5.65,4.35)
\uput[0](3,5){$m$}
\uput[45](1,3){$n$}
\uput[45](1.5,4.5){$\chi_1$}
\uput[135](5.5,4.2){$\chi_2$}
\end{pspicture}
\caption{Incidence \ref{I2} in the case of $m\nparallel n$.}% Right: $m\parallel n$.} 
\label{s2}
\end{figure}

\noindent In this incidence, the reflection of a given line $m$ coincides with another given line $n$. Two cases are possible:

\begin{enumerate}

\item When $m\nparallel n$, there are two possible fold lines that satisfy the incidence, which are the bisectors to the angles defined by $m$ and $n$ (Fig. \ref{s2}).  
\item When $m\parallel n$, there is only one solution, which is a fold line parallel and equidistant to both $m$ and $n$ (Fig. \ref{linereflx}, left, with $m'=n$).
\end{enumerate}

\pagebreak
\begin{incidence}
$\mathcal{F}(P)\in m$, with $P\notin m$.
\label{I3}
\end{incidence}

\noindent In this incidence, the reflection of a given point $P$ is on a given line $m$, and the case in which $P$ is already on $m$ is excluded. It has been shown that the fold lines that satisfy the incidence are  tangents to a parabola with focus $P$ and directrix $m$ \citep{Alperin2000,Martin1998}. 

\begin{figure}[!htb]
\centering

\begin{pspicture}(-3,-3)(4,4)
%\showgrid
\psset{xunit=1.2cm,yunit=1.2cm}
\psaxes[linecolor=gray,labels=none,linewidth=1pt]{->}(0,0)(-1.5,-1.5)(3.5,3)[$x$,0][$y$,0]
\psplot[linewidth=1pt,linecolor=vector,plotpoints=100]{-1.5}{3.2}{x x mul 4 div}
\uput[135]{*0}(0,1){$P(0, 1)$}
\psline(-1.5,-1)(3,-1)
\uput[90]{*0}(-1.5,-1){$m$}
\uput[90]{*0}(-1.5,.6){$\Psi$}
\uput[-90]{*0}(2.5,-1){$P'(t, -1)$}
\psline[linecolor=gray,linewidth=1pt, linestyle=dashed](2.5,-1)(2.5,1.56)
\psline[linecolor=gray,linewidth=1pt, linestyle=dashed](0,1)(2.5,1.56)
\psline[linecolor=gray,linewidth=1pt, linestyle=dashed](0,1)(2.5,-1)
\uput[90]{*0}(2.5,1.56){$T$}
\psline[linecolor=foldline,linewidth=1pt](-.2,-1.82)(3.2,2.44)
\uput[90]{*0}(.5,-.8){$\chi$}
\qdisk(0,1){2pt}
\qdisk(2.5,-1){2pt}
\qdisk(2.5,1.56){2pt}
\qdisk(1.25,0){2pt}
\end{pspicture}
\caption{Incidence \ref{I3}. The fold line $\chi$ is tangent to a parabola with focus $P$ and directrix $m$.} 
\label{parab}
\end{figure}
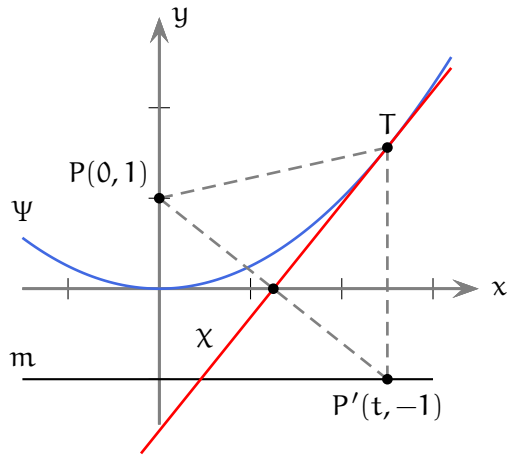

It is useful to derive equations of the fold lines and the associated parabola, for later application to the analysis of fold operations. Without loss of generality, choose a Cartesian system of coordinates $x, y$ so that $P$ is located at $(0, 1)$ and line $m$ is $y=-1$ (Fig.~\ref{parab}). Also, let $P'=\mathcal{F}(P)$ be located at $(t, -1)$, where $t$ is a free parameter. 

The slope of segment $\overline{PP'}$ is $-2/t$. The fold line $\chi$ is perpendicular to $\overline{PP'}$ and therefore has a slope of $t/2$. Further, $\chi$ passes through the midpoint of $\overline{PP'}$, which is located at $(t/2, 0)$. Thus, $\chi$ has an equation 
\begin{equation}
y=\frac{t}{2}\left(x - \frac{t}{2}\right).
\label{chipar}
\end{equation}

Next, consider point $T$ located at the intersection of $\chi$ with a vertical line through $P'$. Its coordinates may be obtained by evaluating Eq.~(\ref{chipar}) at $x=t$, which produces $(t, t^2/4)$. Those coordinates describe parametrically a parabola with equation 
\begin{equation}
y=\frac{x^2}{4},
\label{p1}
\end{equation}
which is denoted by $\Psi$. This is precisely the equation of a parabola with focus at $(0, 1)$ and directrix $y=-1$.\footnote{The general equation of a parabola with vertical axis and vertex at $(0, 0)$ is $y=x^2/(4a)$, where $a$ is the distance from the vertex to the directrix $y=-a$ or the focus $(0, a)$ \citep{Weisstein2016a}.} Further, note that the slope of a tangent to $\Psi$ at point $T$ is $y'(t)=t/2$, which is the same slope of $\chi$. Therefore, $\chi$ is a line tangent to $\Psi$ at point $T$.

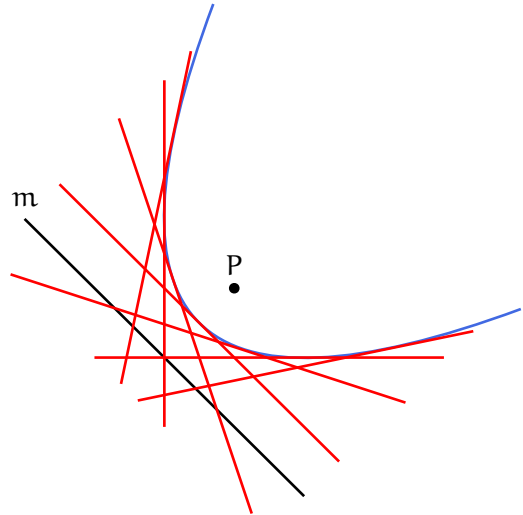
\begin{figure}[!htb]
\centering
\begin{pspicture}(-2,-3)(4,4.5)
%\showgrid
\psset{xunit=1.3cm,yunit=1.3cm}
\rput{-45}{%
\parabola[linecolor=vector,linewidth=1pt](2.2,2.4)(0,0)
\qdisk(0,.5){2pt}
\uput[90]{*0}(0,.5){$P$}
\psline[linewidth=1pt](-2,-.5)(2,-.5)
\uput[90]{*0}(-2,-.5){$m$}
\psline[linecolor=foldline,linewidth=1pt](-2.,0)(2,0)
\psline[linecolor=foldline,linewidth=1pt](-2.02,1.9)(-.12,-1)
\psline[linecolor=foldline,linewidth=1pt](-2,1.5)(.5,-1)
\psline[linecolor=foldline,linewidth=1pt](-2.05,.9)(1.75,-1)
\psline[linecolor=foldline,linewidth=1pt](-1.7,-1)(2.05,.9)
\psline[linecolor=foldline,linewidth=1pt](-.5,-1)(2,1.5)
\psline[linecolor=foldline,linewidth=1pt](.12,-1)(2.02,1.9)
}
\end{pspicture}
\caption{Fold lines for incidence \ref{I3}.}
\label{s3}
\end{figure}

Since $t$ in Eq.~(\ref{chipar}) is a free parameter, then the solution to this incidence is a family of fold lines with one parameter (Fig.~\ref{s3}).

The case $P\in m$ is excluded because under such condition any fold line passing through $P$ or perpendicular to $m$ satisfies the incidence. Those two possibilities are considered in incidences \ref{I4} and \ref{I5}, respectively.

\subsection{Incidences involving an object and its reflected image}

\begin{incidence}
$\mathcal{F}(P)=P$.
\label{I4}
\end{incidence}

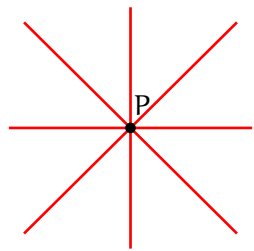
\begin{figure}[!htb]
\centering

\begin{pspicture}(6,1)(10,5)
%\showgrid
\psline[linecolor=foldline,linewidth=1pt](6.6,1.6)(9.4,4.4)
\psline[linecolor=foldline,linewidth=1pt](9.4,1.6)(6.6,4.4)
\psline[linecolor=foldline,linewidth=1pt](8,1.4)(8,4.6)
\psline[linecolor=foldline,linewidth=1pt](6.4,3)(9.6,3)
\qdisk(8,3){2pt}
\uput[80](8.1,3){$P$}
\end{pspicture}
 
\caption{Fold lines for incidence \ref{I4}. } 
\label{s4}
\end{figure}

\noindent In this incidence, the reflection of a given point $P$ coincides with itself, and it is satisfied by any fold line $\chi$ passing through $P$ . An arbitrary direction for line $\chi$ may be defined by its angle $\theta$ with, e.g., the $x$-axis in a Cartesian coordinate system. Therefore, the solution to the incidence is a family of fold lines with one parameter (Fig.~\ref{s4}).

\begin{incidence}
$\mathcal{F}(m)=m$, and $\exists P\in m$, $\mathcal{F}(P)\neq P$.
\label{I5}
\end{incidence}

\noindent Both this and the next incidence consider the reflection of a line $m$ to itself. As discussed in Section \ref{reflections}, there are two ways in which such a reflection may be achieved. In the current case, one half of $m$, defined from an arbitrary point $R\in m$, is reflected upon the opposite half. 

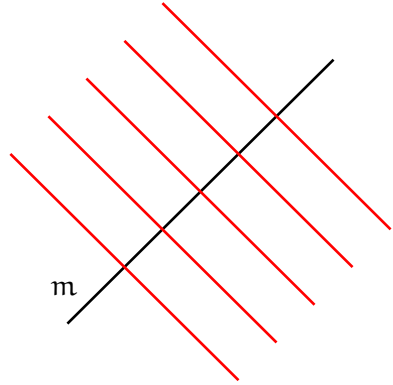
\begin{figure}[!htb]
\centering

\begin{pspicture}(1,0.5)(5,5.5)
%\showgrid
\psline[linewidth=1pt](1.25,1.25)(4.75,4.75)
\psline[linecolor=foldline,linewidth=1pt](2.5,5.5)(5.5,2.5)
\psline[linecolor=foldline,linewidth=1pt](2.,5.)(5,2)
\psline[linecolor=foldline,linewidth=1pt](1.5,4.5)(4.5,1.5)
\psline[linecolor=foldline,linewidth=1pt](1,4)(4,1)
\psline[linecolor=foldline,linewidth=1pt](.5,3.5)(3.5,.5)
\uput[135](1.5,1.5){$m$}
\end{pspicture}
\caption{Fold lines for incidence \ref{I5}.} 
\label{s5}
\end{figure}

The position of point $R$ may be specified by its distance $s$ from a particular point $P_0\in m$. Therefore, the solution to the incidence is a family of fold lines with one free parameter (Fig.~\ref{s5}). 

\begin{incidence}
$\mathcal{F}(m)=m$, and $\forall P\in m$, $\mathcal{F}(P)= P$.
\label{I6}
\end{incidence}

\noindent This is the second case of reflection of line $m$ to itself. In this case, each point $P\in m$ is reflected to itself, and therefore the incidence is satisfied by the unique fold line $\chi=m$.

\section{Elementary fold operations}

\subsection{Definition}
 
A straight line on a plane is an object with two degrees of freedom.\footnote{A fold line $\chi$ may be defined by an equation of the form $ax+by+c=0$, where $a$, $b$ and $c$ are constants, and $(a, b)$ is a normal vector to $\chi$. A vector in arbitrary direction may be defined by  letting $a=\cos \theta$, $b=\sin\theta$, with $0\le \theta < 2\pi$. Therefore, two parameters must be set in order to define any fold line, namely, $\theta$ and $c$.} 
When an incidence constraint is set for $\chi$, satisfying the constraint consumes a number of degrees of freedom, and that number is called the codimension of the constraint. Incidences \ref{I1}, \ref{I2}  and \ref{I6} have either a unique or a finite number of solutions; therefore, those incidences have codimension 2. On 
the other hand, each of incidences  \ref{I3}, \ref{I4} and \ref{I5} have a family of solutions with one free parameter and therefore they have codimension 1 (Table \ref{table0}). 
 
\begin{table}[!htb]
\centering

\begin{threeparttable}
\caption{Incidence constraints.}
\label{table0}
\centering

\begin{tabular}{clc}
\toprule
Incidence & Definition\tnote{a} & Codimension\\
\midrule
\ref{I1} & $\mathcal{F}(P)=Q$, with $P\neq Q$ & 2\\
\ref{I2} & $\mathcal{F}(m)=n$, with $m\neq n$ & 2 \\
\ref{I3} & $\mathcal{F}(P)\in m$, with $P\notin m$ & 1 \\
\ref{I4} & $\mathcal{F}(P)=P$ & 1\\
\ref{I5} & $\mathcal{F}(m)=m$, and $\exists P\in m$, $\mathcal{F}(P)\neq P$ &1\\
\ref{I6} & $\mathcal{F}(m)=m$, and $\forall P\in m$, $\mathcal{F}(P)= P$ & 2\\
\bottomrule
\end{tabular}

\begin{tablenotes}
 \small
\item [a]$P$ and $Q$ are points; $m$ and $n$ are lines.
\end{tablenotes}

\end{threeparttable}

\end{table}

An elementary single-fold operation is defined as a minimal set of alignments between points and lines that is satisfied with a single fold and has a finite number of solutions \citep{Alperin2006}. Equivalently, it is the resolution of a set of incidence constraints which have a total codimension of 2.

Each of the incidences \ref{I1}, \ref{I2} and \ref{I6} already define an elementary operation. The other three incidences must be applied in pairs (including pairing incidences  of the same type), and there are a total of 6 possible pairs. However, incidence \ref{I5} can be not be used twice. If it is, then for given lines $m$ and $n$ the fold line $\chi$ has to satisfy both $\mathcal{F}(m)=m$ and  $\mathcal{F}(n)=n$. Therefore, $\chi$  must be perpendicular to both $m$ and $n$, and two cases are possible:
\begin{enumerate}
\item If $m\nparallel n$, then a perpendicular to both lines does not exist (in the Euclidean plane).  
\item If $m\parallel n$, then any perpendicular to $m$ or $n$ is a valid fold line. 
\end{enumerate}

Thus, the pair of constraints has none or infinite solutions, and so it does not define a valid elementary fold operation.

A total of eight elementary fold operations may be then defined, and they are analyzed in the next subsections.

\subsection{Elementary operations defined by codimension 2 incidences}

The following three operations are defined by incidences \ref{I1}, \ref{I2} and \ref{I6}, respectively. 

\begin{operation}
Given points $P$ and $Q$, with $P\neq Q$, construct a fold line  so that $\mathcal{F}(P)=Q$.
\label{O1}
\end{operation}

\begin{operation}
Given lines $m$ and $n$, with $m\neq n$, construct a fold line so that $\mathcal{F}(m)=n$.
\label{O2}
\end{operation}

\begin{operation}
Given a line $m$, construct a fold line so that $\mathcal{F}(m)=m$, and $\forall P\in m,$ $\mathcal{F}(P)= P$.
\label{O3}
\end{operation}

\subsection{Elementary operations defined by pairs of codimension 1 incidences}

\begin{operation}
Given points $P$ and $Q$, with $P\neq Q$, construct a fold line so that $\mathcal{F}(P)=P$ and $\mathcal{F}(Q)=Q$.
\label{O4}
\end{operation}

\begin{figure}[!htb]
\centering

\begin{pspicture}(1.5,1.5)(5,5.5)
%\showgrid
\uput[135](2.3,2.3){$P$}
\uput[135](4.5,4.5){$Q$}
\psline[linecolor=foldline,linewidth=1pt](1.8,1.8)(5,5)
\qdisk(2.3,2.3){2pt}
\qdisk(4.5,4.5){2pt}
\uput[135](3.5,3.5){$\chi$}
\end{pspicture}
\caption{Operation \ref{O4}.}
\label{fO4}
\end{figure}
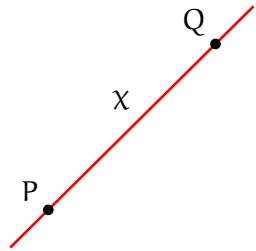

\noindent
This operation is defined by application of incidence \ref{I4} to two distinct points. It has a unique solution, which is a fold line $\chi$ passing through points $P$ and $Q$ (Fig.~\ref{fO4}).

\begin{operation}
Given a point $P$ and a line $m$, construct a fold line such that $\mathcal{F}(P)=P$ and $\mathcal{F}(m)=m$, and $\exists Q\in m$, $\mathcal{F}(Q)\neq Q$.
\label{O5}
\end{operation}

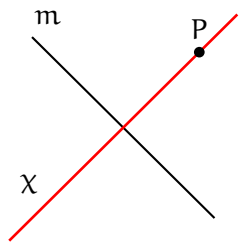
\begin{figure}[!htb]
\centering

\begin{pspicture}(-2,-2)(2,2.5)
%\showgrid
\psline(-1.2,1.2)(1.2,-1.2)
\uput[90](-1.,1.2){$m$}
\psline[linecolor=foldline,linewidth=1pt](-1.5,-1.5)(1.5,1.5)
\qdisk(1,1){2pt}
\uput[90]{*0}(1,1){$P$}
\uput[135]{*0}(-1,-1){$\chi$}
\end{pspicture}
\caption{Operation \ref{O5}.}
\label{fO5} 
\end{figure}

\noindent
This operation is defined by application of incidences \ref{I4} and \ref{I5}. It has a unique solution, which is a fold line $\chi$ perpendicular to $m$ and passing through $P$ (Fig.~\ref{fO5}). Note that the case $P\in m$ is allowed, which has the same unique solution.

\begin{operation}
Given points $P$, $Q$, and a line $m$, with $P\notin m$, construct a fold line such that $\mathcal{F}(P)\in m$ and $\mathcal{F}(Q)=Q$.
\label{O6}
\end{operation}

\noindent 
This operation is defined by application of incidences \ref{I3} and \ref{I4}. Its solution is a fold line that is tangent to a parabola with focus $P$ and directrix $m$, and passes through point $Q$. 

Assume the same parabola of Fig.~\ref{parab}, given by Eq.~(\ref{p1}), and a point $Q$ at the position $(x_q, y_q)$. Replacing the coordinates of $Q$ in Eq.~(\ref{chipar}) produces the quadratic equation
\begin{equation}
t^2-2x_qt+4y_q=0.
\label{tyt}
\end{equation}
The discriminant of Eq.~(\ref{tyt}) is $\Delta=4x_q^2-16y_q$, and $\Delta=0$ yields $y_q=x_q^2/4$, which implies $Q\in \Psi$. Since $\Psi$ is the location of points that are equidistant from $P$ and  $m$, we may conclude that the fold operation has a unique solution when $Q$  is equidistant to $P$ and $m$, two solutions when $Q$ is closer to $m$ (i.e., $y_q<x_q^2/4$), and no solution when $Q$ is closer to $P$ (i.e., $y_q>x_q^2/4$).

\begin{figure}[!htb]
\centering

\begin{pspicture}(-3,-2.5)(3,4.5)
%\showgrid
\rput{-45}{%
\psplot[linecolor=vector,plotpoints=100,linewidth=1pt]{-2.5}{2.5}{x x mul 2 div}
\qdisk(0,.5){2pt}
\uput[90]{*0}(0,.5){$P$}
\psline(-2.5,-.5)(2.5,-.5)
\uput[-135]{*0}(-2,-.5){$m$}
\psline[linecolor=foldline,linewidth=1pt](-2.5,.53)(2.5,-.6)
\psline[linecolor=foldline,linewidth=1pt](-.4,-1.5)(2.5,2.45)
\qdisk(0.57,-.17){2pt}
\uput[-30]{*0}(0.7,-.17){$Q$}
\uput[0]{*0}(-1.8,1.56){$\Psi$}
\uput[-135]{*0}(-2,.6){$\chi_1$}
\uput[-45]{*0}(-.2,-1.3){$\chi_2$}
}
\end{pspicture}
\caption{Example of operation \ref{O6} with two solutions.}
\label{fO6} 
\end{figure}
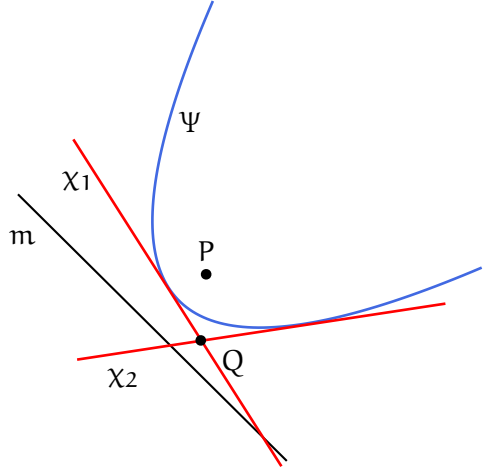

Fig.~\ref{fO6} shows an example for the case of two solutions.

\begin{operation}Given points $P$, $Q$, and lines $m$, $n$, with $P\notin m$, $Q\notin n$, and $P\neq Q$ or $m \neq n$, construct a fold line so that $\mathcal{F}(P)\in m$ and $\mathcal{F}(Q)\in n$.
\label{O7}
\end{operation}

\noindent
This operation derives from the application of incidence \ref{I3} to two distinct point-line pairs. Its solution is a fold line that is tangent to both a parabola $\Psi$ with focus $P$ and directrix $m$, and a parabola $\Theta$ with focus $Q$ and directriz $n$. 

Again, assume the same parabola $\Psi$ of Fig.~\ref{parab}, given by Eq.~(\ref{p1}). Assume also that  $Q$ is located at $(x_q, y_q)$, and its reflection $Q'=\mathcal{F}(Q)$ is at $(x_q', y_q')$. Then, segment $\overline{QQ'}$ has a slope $(y_q - y_q')/(x_q-x_q')$. 
The fold line $\chi$, given by Eq.~(\ref{chipar}), has a slope $t/2$ and is perpendicular to $\overline{QQ'}$ (because it reflects $Q$ onto $Q'$). Therefore, 
\begin{equation}
\frac{t}{2}=-\frac{x_q-x_q'}{y_q-y_q'}.
\label{o61}
\end{equation}
Further, $\chi$ passes through the midpoint of $\overline{QQ'}$, which is located at $((x_q+x_q')/2, (y_q + y_q')/2)$. Replacing those coordinates into Eq.~(\ref{chipar}) produces  
\begin{equation}
2(y_q+y_q') = t\left(x_q+x_q'-t\right).
\label{o62}
\end{equation}
Finally, eliminating $t$ from Eqs.~(\ref{o61}) and (\ref{o62}) produces
\begin{equation}
(y_q+y_q')(y_q-y_q')^2=-(x_q^2-x_q'^2)(y_q-y_q')-2(x_q-x_q')^2.
\label{cubic}
\end{equation}

For a given line $n$, the coordinates of $Q'$ satisfy an equation of the form
\begin{equation}
ax_q'+by_q'+c=0,
\label{linen}
\end{equation}
where $a$, $b$ and $c$ are constants.  

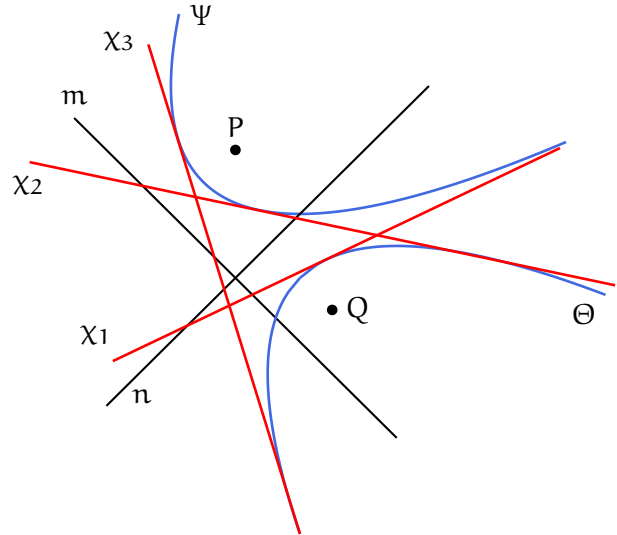
\begin{figure}[!htb]
\centering

\psset{xunit=1.2cm,yunit=1.2cm}
\begin{pspicture}(-2,-4)(4,2.5)
%\showgrid
\rput{-45}{%
\psplot[linecolor=vector,plotpoints=100,linewidth=1pt]{-1.5}{2.5}{x x mul 2 div}
\qdisk(0,.5){2pt}
\uput[90]{*0}(0,.5){$P$}
\psline(-1.5,-.5)(3.5,-.5)
\uput[90]{*0}(-1.5,-.5){$m$}
\qdisk(2,0){2pt}
\uput[0]{*0}(2,0){$Q$}
\psline(1,-2.5)(1,2.5)
\uput[-90]{*0}(1.05,-2){$n$}
\psplot[linecolor=vector,plotpoints=100,linewidth=1pt]{1.5}{4}{x 2 mul 3 sub sqrt}
\psplot[linecolor=vector,plotpoints=100,linewidth=1pt]{1.5}{3.5}{x 2 mul 3 sub sqrt -1 mul}
\psplot[linecolor=foldline,plotpoints=100,linewidth=1pt]{-1.5}{3.5}{x -1.5 add -.53 mul -.94 add}
\psplot[linecolor=foldline,plotpoints=100,linewidth=1pt]{-1.5}{4}{x -1.5 add .65 mul 0.76 add}
\psplot[linecolor=foldline,plotpoints=100,linewidth=1pt]{.7}{2.5}{x -1.5 add 2.85 mul 0.18 add}
\uput[-90]{*0}(-1.5,-1.2){$\chi_2$}
\uput[135]{*0}(-1.4,0.5){$\chi_3$}
\uput[180]{*0}(0.6,-1.8){$\chi_1$}
\uput[0]{*0}(-1.5,1.1){$\Psi$}
\uput[-90]{*0}(3.8,2.1){$\Theta$}
}
\end{pspicture}
\caption{Example of operation \ref{O7} with three solutions.} 
\label{fO7}
\end{figure}

Eqs.~(\ref{cubic}) and (\ref{linen}) may be solved for $x_q'$ and $y_q'$. Substituting this solution into Eq.~(\ref{o61}) yields $t$, which defines the fold line $\chi$ in Eq.~(\ref{chipar}). Two cases may be considered: 
\begin{enumerate}
\item If $m\parallel n$, then $Q'$ is on a horizontal line and so $y_q' = -c/b$. In that case, Eq.~(\ref{cubic}) is quadratic in $x_q'$ and may have zero to two solutions.

\item If $m\nparallel n$, solving Eq.~(\ref{linen}) for $x_q'$ or $y_q'$ and replacing in Eq.~(\ref{cubic}) produces a cubic equation with one to three solutions. An example for the latter case is shown in Fig.~\ref{fO7}. 
\end{enumerate}

Let us investigate further the conditions to ensure the existence of solutions. As noted above, the operation may not have a solution only in the case of $m\parallel n$. Re-arranging Eq.~(\ref{cubic}) produces
\begin{equation}
(2-y_q+y_q')(x_q-x_q')^2+2x_q(y_q-y_q')(x_q-x_q')+(y_q+y_q')(y_q-y_q')^2=0,
\label{cubic2}
\end{equation}
which is a quadratic equation in $(x_q-x_q')$. The discriminant is 
\begin{equation}
\Delta=4x_q^2(y_q-y_q')^2-4(2-y_q+y_q')(y_q+y_q')(y_q-y_q')^2,
\end{equation}
and letting $\Delta\ge 0$ produces
\begin{equation}
x_q^2+(y_q-1)^2\ge (y_q'+1)^2.
\label{cubic3}
\end{equation}

The left side of Eq.~(\ref{cubic3}) is the squared distance between $P$ and $Q$, and the right side is the squared distance between $m$ and $n$.  This result does not seem reported in the literature, and may be stated as a theorem:

\begin{theorem}
Given points $P$, $Q$, and lines $m$, $n$, with $P\notin m$, $Q\notin n$, and $P\neq Q$ or $m \neq n$, a fold line that places $P$ on $m$ and $Q$ on $n$ exists iff the distance between $P$ and $Q$ is larger than or equal to the distance between $m$ and $n$.
\end{theorem}

\begin{operation}
Given point $P$ and lines $m$ and $n$, with $P\notin m$, construct a fold line so that $\mathcal{F}(P)\in m$ and $\mathcal{F}(n)=n$ and $\exists Q\in n, \mathcal{F}(Q)\neq Q$.
\label{O8}
\end{operation}

\noindent 
This operation is defined by application of incidences \ref{I3} and \ref{I5}. Its solution is a fold line that is tangent to a parabola with focus $P$ and directrix $m$, and is perpendicular to line $n$. 

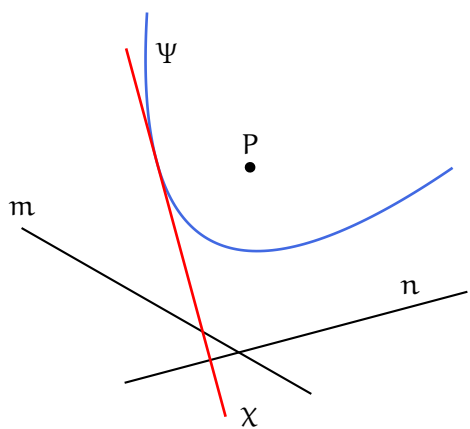
\begin{figure}[!htb]
\centering

\begin{pspicture}(-3,-2.5)(3,3)
%\showgrid
\psset{xunit=1.1cm,yunit=1.1cm}
\rput{-30}{%
\psplot[linecolor=vector,linewidth=1pt,plotpoints=100]{-2}{2.1}{x x mul 2 div}
\qdisk(0,1){2pt}
\uput[90]{*0}(0,1){$P$}
\psline(-2,-1)(2,-1)
\uput[90]{*0}(-2,-1){$m$}
\psplot[plotpoints=10]{0}{3}{x -2 add}
\psplot[linecolor=foldline,linewidth=1pt,plotpoints=10]{-2}{1.25}{x -1 mul -.5 add}
\uput[90]{*0}(2.5,.5){$n$}
\uput[0]{*0}(-1.8,1.56){$\Psi$}
\uput[-30]{*0}(1.2,-1.6){$\chi$}
}
\end{pspicture}
\caption{Example of operation \ref{O8} with one solution.}
\label{fO8} 
\end{figure}

As in the previous operation, assume the same parabola of Fig.~\ref{parab} given by Eq.~(\ref{p1}), and a line $n$ with equation $ax+by+c=0$. The fold line $\chi$ has a slope $t/2$ and is perpendicular to $n$. Two cases may be considered:
\begin{enumerate}
\item If $m\nparallel n$, then $a\neq 0$. Therefore,
\begin{equation}
\frac{t}{2}=-\frac{b}{a}
\end{equation}
which has a unique solution for $t$ (Fig.~\ref{fO8}). Knowing $t$, Eq.~(\ref{chipar}) defines the fold line $\chi$. 

\item If $m\parallel n$ then $n$ is a horizontal line and cannot be perpendicular to any tangent to parabola $\Psi$. In this case, the operation does not have a solution.
\end{enumerate}

\subsection{Summary}

Table \ref{table2} lists the incidence constraints that define each operation and their number of solutions, Table \ref{table2b} lists the conditions for the existence of solutions, and Table \ref{table3} restates the operations as folding actions of a medium $\mathcal{O}$.

\begin{table}[!htb]
\centering
\begin{threeparttable}
\caption{Incidence constraints and number of solutions of the elementary single-fold operations.}
\label{table2}
\centering
\begin{tabularx}{\textwidth}{lXl}
\toprule
Operation & Incidence constraints & Solutions\\
\midrule
\ref{O1} & $\mathcal{F}(P)=Q$, with $P\neq Q$ & 1\\
\ref{O2} & $\mathcal{F}(m)=n$, with $m\neq n$ & 1, 2\\
\ref{O3} & $\mathcal{F}(m)=m$,and $\forall P\in m, \mathcal{F}(P)=P$  &1\\
\ref{O4} & $\mathcal{F}(P)=P$ and $\mathcal{F}(Q)=Q$, with $P\neq Q$ & 1\\
\ref{O5} & $\mathcal{F}(P)=P$ and $\mathcal{F}(m)=m$, and $\exists Q\in m, \mathcal{F}(Q)\neq Q$ & 1\\
\ref{O6} & $\mathcal{F}(P)\in m$, with $P\notin m$, and $\mathcal{F}(Q)=Q$ & 0 -- 2\\
\ref{O7} & $\mathcal{F}(P)\in m$, with $P\notin m$, $\mathcal{F}(Q)\in n$, with $Q\notin n$, and $P\neq Q$ or $m\neq n$ & 0 -- 3\\
\ref{O8} &  $\mathcal{F}(P)\in m$, with $P\notin m$, and $\mathcal{F}(n)=n$, and $\exists Q\in n$, $\mathcal{F}(Q)\neq~Q$ & 0, 1\\
\bottomrule
\end{tabularx}
\end{threeparttable}
\end{table}

\begin{table}[!htb]
\centering
\begin{threeparttable}
\caption{Conditions for the existence of solutions of the elementary single-fold operations.}
\label{table2b}
\centering
\begin{tabularx}{\textwidth}{lX}
\toprule
Operation & Conditions\\
\midrule
\ref{O1} to \ref{O5} & none\\
\ref{O6} & distance between $P$ and $Q$ larger than or equal to the distance between $Q$ and $m$\\
\ref{O7} & distance between $P$ and $Q$ larger than or equal to the distance between $m$ and $n$\\
\ref{O8} &  $m\nparallel n$\\
\bottomrule
\end{tabularx}
\end{threeparttable}
\end{table}

\begin{table}[!htb]
\centering
\begin{threeparttable}
\caption{Elementary single-fold operations restated as folding actions.}
\label{table3}
\centering
\begin{tabularx}{\textwidth}{cX}
\toprule
No. &\multicolumn{1}{c}{Action\tnote{a}}\\
\midrule
\ref{O1} & Given two distinct points $P$ and $Q$, fold $\mathcal{O}$ to place $P$ onto $Q$.\\
\ref{O2} &  Given two distinct lines $m$ and $n$, fold $\mathcal{O}$ to align $m$ and $n$.\\
\ref{O3} & Fold along a given a line $m$.\\
\ref{O4} & Given two distinct points $P$ and $Q$, fold $\mathcal{O}$ along a line passing through $P$ and $Q$.\\
\ref{O5} & Given a line $m$ and a point $P$, fold $\mathcal{O}$ along a line passing through $P$ to reflect half of $m$ onto its other half.\\
\ref{O6} & Given a line $m$, a point $P$ not on $m$ and a point $Q$, fold $\mathcal{O}$ along a line passing through $Q$ to place $P$ onto $m$.\\
\ref{O7} & Given two lines $m$ and $n$, a point $P$ not on $m$ and a point $Q$ not on $n$, where $m$ and $n$ are distinct or $P$ and $Q$ are distinct, fold $\mathcal{O}$ to place $P$ onto $m$, and $Q$ onto $n$.\\
\ref{O8} & Given two lines $m$ and $n$, and a point $P$ not on $m$, fold $\mathcal{O}$ to place $P$ onto $m$, and to reflect half of $n$ onto its other half.\\
\bottomrule
\end{tabularx}
\begin{tablenotes}
 \small
\item [a]$\mathcal{O}$ denotes the medium in which folds are performed; e.g., a sheet of paper, fabric, plastic, metal or any other foldable material.
\end{tablenotes}
\end{threeparttable}
\end{table}

\section{Discussion}
\label{foldingalong}

The complete set of elementary single-fold operations contains eight operations, listed in Table \ref{table3}. Operations \ref{O1}, \ref{O2} and \ref{O4} to \ref{O8} constitute Justin's  original set \citep{Justin1986}, and operation \ref{O3} is the new addition proposed here. 
\ref{O3} does not create a new line and has been ignored in previous studies on origami constructions \citep{Alperin2006, Ghourabi2013}.\footnote{In his formulation, \citet{Justin1986} allowed for solutions where the fold line itself coincides with an existent line. However, such action was not considered as a fold operation on its own.} However, it is a valid  elementary single-fold operation and completeness of the set demands its inclusion.

There is also a more practical reason for not ignoring \ref{O3}. Folding a sheet of paper along a line superposes the paper on both sides of the fold line, in two layers. In origami mathematics, it is assumed that all lines and points marked on one layer are also defined on the layers above and below, as if the paper were ``transparent'' \citep{Geretschlager1995, Martin1998}. However, it is not so in actual paper folding: a given origami work may require to fold, e.g., the top layer along a line marked on the layer below it. Such is the case when folding parallel lines for building tessellation grids \citep[see instruction 5 for a triangle grid in p. 8]{Gjerde2008}. 

\begin{figure}[!htb]
\centering
\begin{pspicture}(0,-1)(15,4)
\psset{xunit=.8cm,yunit=.8cm}
%\showgrid
\psframe[linewidth=0.5pt](0,0)(4,4) 
\uput[45](.5,0){$m$}
\uput[45](1.5,0){$n$}
\psline[linewidth=1pt](.5,0.1)(1.5,3.9)
\psline[linewidth=1pt](1.5,0.1)(2.5,3.9)
\psarc[linecolor=darkgray]{->}(1.5,0){2}{40}{110}

\pspolygon[linewidth=0.5pt](6,0)(7.5,0)(8.5,4)(6,4) 
\psline[linewidth=1pt](6.5,0.1)(7.5,3.9)
\pspolygon[linewidth=0.5pt,fillstyle=solid,fillcolor=backside](5.3,1.25)(7.5,0)(8.5,4)(7.25,4.7) 
\psline[linestyle=dashed,linecolor=gray](6.62,0.5)(7.5,3.9)
\psline[linestyle=dashed,linecolor=gray](6.9,4)(8.5,4)
\uput[135](6.6,0){$m$}
\uput[45](7.5,0){$n$}
\psarcn[linecolor=darkgray]{->}(6.7,0.5){1.75}{95}{60}

\pspolygon[linewidth=0.5pt](10,0)(11.5,0)(12.5,4)(10,4) 
\psline[linewidth=1pt](10.53,0.1)(10.62,.45)
\pspolygon[linewidth=0.5pt,fillstyle=solid,fillcolor=backside](10.62,.45)(11.5,0)(12.5,4)(11.62,4.5) 
\pspolygon[linewidth=0.5pt,fillstyle=solid,fillcolor=white](10.61,.45)(12.1,0.45)(12.1,4.45)(11.61,4.5) 

\uput[135](10.6,0){$m$}
\uput[45](11.5,0){$n$}
\psarcn[linecolor=darkgray]{->}(11.7,0.){1.75}{95}{60}

\psframe[linewidth=0.5pt](14,0)(18,4) 
\uput[45](14.5,0){$m$}
\uput[45](15.5,0){$n$}
\uput[45](16.5,0){$o$}
\psline[linewidth=1pt](14.5,0.1)(15.5,3.9)
\psline[linewidth=1pt](15.5,0.1)(16.5,3.9)
\psline[linewidth=1pt](16.5,0.1)(17.5,3.9)

\uput[-90](2,-.3){(1)}
\uput[-90](7.25,-.3){(2)}
\uput[-90](11.25,-.3){(3)}
\uput[-90](16,-.3){(4)}

\end{pspicture}
 
\caption{Given parallel lines $m$ and $n$, create a third equidistant parallel line at the left of $n$. (1) Fold along line $n$. (2) Fold the top layer along line $m$ in the bottom layer. (3) Unfold. (4) Final result.} 
\label{parf}

\end{figure}
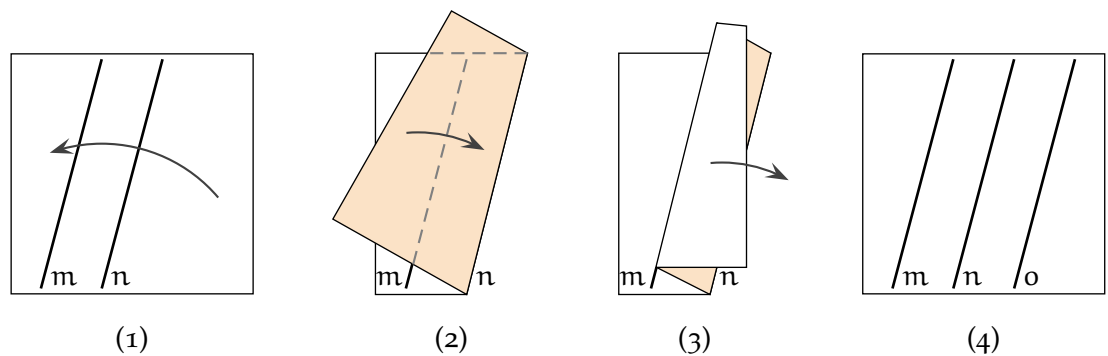

Fig.~\ref{parf} shows a simple example. Assume a sheet of paper (not necessarily square) in which two parallel lines $m$ and $n$ are marked, and assume that we want to create a third equidistant parallel line $o$. Both steps (1) and (2) require to fold along given lines ($n$ and $m$, respectively). Naturally, it would be also possible to use the points of intersections of $m$ and $n$ with the borders of the paper as references, instead of the lines themselves. Thus, the instruction for step (2) could be: fold the top layer along a line passing through points $P$ and $Q$ on the bottom layer, where $P$ and $Q$ are the intersections of $m$ with the upper and lower edges of the paper. However,  in actual practice it is much simpler and convenient to perform the fold by aligning it with line $m$. Further, it might be the case that the borders of the paper are not well defined (or have not been defined) or that it is considered as a theoretical infinite plane. 

\ref{O3} is also a common instruction also for building figurative models \citep[see steps 8 and 16 of ``Baby'' in page 88, and step 6 of ``Songbird 2'' in page 340]{Lang2012}.  Thus, the operation should be in the repertoire of, e.g., computational systems for origami simulation and design \citep{Ida2009,Tsuruta2009}.   

\section{Conclusion}

Analysis of reflections of points and lines on a plane subject to incidence constraints has determined a complete set of eight elementary single-fold operations. Precise definitions and conditions of existence of solutions of all operations are given in Tables \ref{table2} to \ref{table3}, which may be useful to scientific and technological applications of origami. 
 
%\bibliographystyle{abbrv}
%\bibliography{../../origami}

\end{document}